\documentclass{article}
\usepackage{amssymb}
\newcommand{\bC}{\mathbf{C}}
\newcommand{\bN}{\mathbf{N}}
\newcommand{\bR}{\mathbf{R}}

\newcommand{\ord}{\mbox{\rm ord}}
\newcommand{\Nj}{{\cal N}_J}

\newcommand{\bb}{\bar{b}}
\newcommand{\jac}{\mathrm{jac}\,}
\newcommand{\supp}{\mathrm{supp}\,}
\newcommand{\Discr}{\mathrm{Discr}_{y}}

\newtheorem{Theorem}{Theorem}[section]
\newtheorem{Remark}[Theorem]{Remark}
\newtheorem{Example}[Theorem]{Example}
\newtheorem{Lemma}[Theorem]{Lemma}
\newtheorem{Formula}[Theorem]{Formula}
\newtheorem{Corollary}[Theorem]{Corollary}

\newenvironment{proof}[1][Proof]{\textbf{#1.} }{\
\rule{0.5em}{0.5em}}

\newcommand{\Teis}[2]{\left\{
   \setlength{\unitlength}{1ex}
   \begin{picture}(2,0)(0,0.4)
      \put(0,1.1){\line(1,0){2}}
      \put(0,0.9){\line(1,0){2}}
      \put(1,1.2){\makebox(0,0)[b]{$\scriptstyle #1$}}
      \put(1,0.8){\makebox(0,0)[t]{$\scriptstyle #2$}}
   \end{picture}\right\}}

\newcommand{\Teisssr}[4]{\left\{
   \setlength{\unitlength}{1ex}
   \begin{picture}(#3,3)(0,0.4)
      \put(0,1.15){\line(1,0){#3}}
      \put(0,0.85){\line(1,0){#3}}
      \put(#4,1.3){\makebox(0,0)[b]{$#1$}}
      \put(#4,0.7){\makebox(0,0)[t]{$#2$}}
   \end{picture}\right\}}

\begin{document}
\title{A discriminant criterion of irreducibility
\footnotetext{
     \noindent   \begin{minipage}[t]{5in}
       {\small
       2000 {\it Mathematics Subject Classification:\/} Primary 32S55;
       Secondary 14H20.\\
       Key words and phrases: irreducible plane curve,
       jacobian Newton polygon, discriminant curve, branches at infinity.\\
       The first-named author was partially supported by the Spanish Project
       PNMTM 2007-64007.}
       \end{minipage}}}
       
\author{Evelia R.\ Garc\'{\i}a Barroso and Janusz Gwo\'zdziewicz}

\maketitle

\begin{abstract}
\noindent In this paper we give a criterion of irreducibility for a
complex power series in two variables, using the notion of jacobian
Newton diagrams, defined with respect to any direction. Moreover we
study the singularity at infinity of a plane affine curve with one
point at infinity for which the global counterpart of our main
result holds.
\end{abstract}

\section{Introduction}
\noindent In \cite{GG} we give criteria of irreducibility for a
complex power series in two variables, using the notion of jacobian
Newton diagrams, defined with respect to a generic direction. In
this paper we generalize these criteria to any direction and we use
this new general criterion to study the branches at infinity of
polynomial curves. The paper is organized as follows:

\medskip

\noindent In~\ref{S:Nd} we recall the notion of the {\em Newton
diagram.} Then in~\ref{S:Dc} we explain what we mean by the {\em
discriminant curve} of an analytic mapping
$F:(\bC^2,0)\to(\bC^2,0)$. If $D(u,v)=0$ is an equation of the
discriminant curve then the Newton diagram of $D$ will be called
{\em jacobian Newton diagram} of $F$ and denoted~$\Nj(F)$. At the
end of the section we present formulas for computing equations of
discriminants.

\medskip

\noindent In Section \ref{jacobian branch} we consider $\Nj(l,f)$
where $l$ is a regular function and $f$ is a singular irreducible
series. We shall call such diagrams the {\em Merle type diagrams}.
We recall Merle's result that equisingularity class of $f$ and the
intersection multiplicity $(f,l)_0$ determine and are determined by
$\Nj(l,f)$. In Theorem~\ref{Th1} we give necessary and sufficient
conditions of arithmetical nature for a Newton diagram to be a Merle
type diagram.

\medskip

\noindent The main result of the paper is Theorem~\ref{Th2}. It
states that $f$ is an irreducible power series if and only if
$\Nj(l,f)$ is a Merle type diagram. We apply our irreducibility
criterion to power series taken from Kuo's paper \cite{Kuo}.

\medskip

\noindent Finally we study the singularity at infinity of a plane
affine curve with one point at infinity for which the global
counterpart of our main result holds. This is interesting in the
context of the Jacobian conjecture. Recall that Abhyankar proved in
\cite{Abhyankar1} that this conjecture is settled affirmatively in
the case where there is only one branch at infinity.

\subsection{Newton diagrams of plane analytic curves}\label{S:Nd}

\noindent In this section we recall the notion of a Newton diagram
and establish the notation. Write $\bR_{+}=\{\,x\in\bR: x\geq0\,\}$.

\medskip

\noindent Let $f\in\bC\{x,y\}$, $f(x,y)=\sum a_{i,j}x^iy^j$ be a
non-zero convergent power series. Put $\supp f:=\{\,(i,j): a_{i,j}
\neq0\,\}$. Then by definition the {\em Newton diagram} $\Delta_f$,
in the coordinates $(x,y)$, of $f$ is
$$
\Delta_f = \mbox{Convex Hull}\;(\supp f+\bR_{+}^2).
$$

\noindent The basic property of Newton diagrams is that the Newton
diagram of a product is the Minkowski sum of Newton diagrams. There
is $\Delta_{fg}=\Delta_f+\Delta_g$ where
$\Delta_f+\Delta_g=\{\,a+b:a\in\Delta_f, b\in\Delta_g\,\}$. In
particular if $f$ and $g$ differ by an invertible factor $u\in
\mathbf C\{x,y\}$, $u(0,0)\neq0$ then $\Delta_f=\Delta_g$. Thus the
Newton diagram of a plane analytic curve is well defined because an
equation of an analytic curve is given up to invertible factor.

\medskip

\noindent After Teissier \cite{Teissier2} we introduce {\em
elementary Newton diagrams}. For $m,n>0$  we put
$\Teis{n}{m}=\Delta_{x^n+y^m}$. We put also
$\Teis{n}{\infty}=\Delta_{x^n}$ and $\Teis{\infty}{m}=\Delta_{y^m}$.

\medskip

\noindent One can check that every Newton diagram $\Delta \subsetneq
\bR_{+}^2$ has a unique representation $\Delta=\sum_{i=1}^r
\Teis{L_i}{M_i}$ where {\em inclinations} of successive elementary
diagrams form an increasing sequence (by definition the inclination
of  $\Teis{L}{M}$ is $L/M$ with the conventions that $L/\infty=0$
and $\infty/M=+\infty$). We shall call this representation the {\em
canonical form} of $\Delta$.

\medskip

\noindent Finally a Newton diagram is {\em convenient} if it
intersects both coordinate axes.

\subsection{Discriminant curve}\label{S:Dc}
\noindent Let $F=(p,q):(\bC^2,0)\to(\bC^2,0)$ be an analytic mapping
given by $(u,v)=(p(x,y),q(x,y))$ and such that
$F^{-1}(0,0)=\{(0,0)\}$. Let $\jac(p,q)=0$ be the equation of the
critical locus of $F$ where $\jac(p,q)=\frac{\partial p}{\partial
x}\frac{\partial q}{\partial y} - \frac{\partial p}{\partial
y}\frac{\partial q}{\partial x}$ is the usual jacobian determinant.
The direct image of~$\jac(p,q)=0$ by $F$ is called the {\em
discriminant curve} of $F$ (see \cite{Casas2}).

\medskip

\noindent Assume that $D(u,v)=0$ is an equation of the discriminant
curve. Then $\Delta_D$, in the coordinates $(u,v)$, is called the
{\em jacobian Newton diagram} of $F$ (see \cite{Teissier3}). We will
write $\Nj(p,q)$ for the jacobian Newton diagram.

\medskip

\noindent Below we give some formulas for jacobian Newton diagrams
and discriminant curves.

\begin{Formula}[Teissier's formula \cite{Teissier1}]
Assume that $\jac(p,q)=h_1\cdots h_r$ where $h_i$ are irreducible
series for $1\leq i\leq r$. Then
$$ \Nj(p,q)=\sum_{i=1}^r\Teisssr{(q,h_i)_0}{(p,h_i)_0}{8}{4}\;
$$

\noindent where $(f,g)_0$ denotes the intersection number of $f$ and
$g$.
\end{Formula}

\noindent From now on we will only consider the mappings
\begin{equation}\label{Eq3}
  (l,f):(\bC^2,0)\to(\bC^2,0)
\end{equation}
where $l$ is a regular function (i.e. $l=ax+by+\mbox{higher order
terms}$, $ax+by\not\equiv0$) and $f$ is a singular series. Recall
that a power series is called \textit{singular} if its order is
larger than one. Under these assumptions $\jac(l,f)=0$ is called the
{\em polar curve of f with respect to l}. The inclinations of the
elementary diagrams of the jacobian Newton diagram $\Nj(l,f)$ are
called {\em polar quotients}. These notions were studied by many
authors (see for example \cite{Merle} and \cite{Ephraim} for
irreducible case, and \cite{Casas}, \cite{Delgado1},
\cite{Delgado2}, \cite{Eggers}, \cite{GB}, \cite{GP3},
\cite{Maugendre1}, \cite{Maugendre2} and \cite{Wall} among others
for the reduced case). Also~\cite{GLP} is a survey of recent
results. If the curves $l=0$ and $f=0$ are transverse then
$\Nj(l,f)$ depends only on the equisingularity class of $f=0$ (see
\cite{Teissier1}). Otherwise the jacobian Newton diagram may depend
on relative position of curves $l=0$ and $f=0$ as the following
example shows.

\begin{Example}\label{Ex:1} Let $f=y^2-x^5$ and let $l_1=x$, $l_2=y$,
$l_3=y-x^2$. Then $\jac(l_1,f)=2y$, $\jac(l_2,f)=5x^4$ and
$\jac(l_3,f)=x(5x^3-4y)$. By Teissier's formula
\begin{eqnarray*}
\Nj(l_1,f) &=&
  \Teisssr{(f,y)_0}{(l_1,y)_0}{8}{4}=\Teisssr{5}{1}{3}{1.5} \\
\Nj(l_2,f) &=&
  4\Teisssr{(f,x)_0}{(l_2,x)_0}{8}{4}=4\Teisssr{2}{1}{3}{1.5}=\Teisssr{8}{4}{3}{1.5}
  \\
\Nj(l_3,f) &=&
\Teisssr{(f,x)_0}{(l_3,x)_0}{8}{4} +
\Teisssr{(f,5x^3-4y)_0}{(l_3,5x^3-4y)_0}{16}{8} =
\Teisssr{2}{1}{3}{1.5} + \Teisssr{5}{2}{3}{1.5} \;.
\end{eqnarray*}
\end{Example}

\noindent For any local analytic diffeomorphism
$\Phi:(\bC^2,0)\to(\bC^2,0)$ the subs\-titution
$(l_1,f_1)=(l\circ\Phi,f\circ\Phi)$ does not affect the equation of
the discriminant curve. Hence without loss of generality we may
assume that $l=x$ (take such a $\Phi$ that $l\circ\Phi=x$).

\begin{Formula}\label{Rem2}
If $f(x,y)$ is a convergent power series such that
$f(0,y)=y^n+\mbox{higher order terms}$ and $\frac{\partial
f}{\partial y}(x,y) = \mathit{unit}\prod_{j=1}^{n-1}[y-\gamma_j(x)]$
is a Newton-Puiseux factorization of $\frac{\partial f}{\partial y}$
then the discriminant of the mapping $(x,f):(\bC^2,0)\to(\bC^2,0)$
has, up to an invertible factor, an equation
\begin{equation}\label{Eq:Rem2}
   D(u,v)=\prod_{j=1}^{n-1}[v-f(u,\gamma_j(u))].
\end{equation}
\end{Formula}

\noindent See the Appendix for the proof.

\begin{Formula}\label{Rem1}
If $f(x,y)=y^n+a_1(x)y^{n-1}+\dots+a_n(x)\in\bC\{x\}[y]$ is a
Weierstrass polynomial, i.e. $a_i(0)=0$ for every $i\in
\{1,\ldots,n\}$, then the discriminant of the mapping
$(x,f):(\bC^2,0)\to(\bC^2,0)$ has an equation $D(u,v)=0$ where
\begin{equation}\label{Eq:Rem1}
   D(u,v)=\Discr(f(u,y)-v)
\end{equation}
is the classical discriminant of a polynomial in one variable~$y$.
\end{Formula}
\begin{proof}
The discriminant $\Discr(f(u,y)-v)$ is, up to an integer constant,
equal to the resultant of polynomials~$f(u,y)-v$ and $\frac{\partial
f}{\partial y}(u,y)$. By the classical formula (see Theorem 10.10,
Chapter I,~\cite{Walker}) the resultant of polynomials $P,Q\in
K[Y]$, $Q=\prod_{i=1}^s(Y-\beta_i)$ where $K$ is a field is, up to a
sign, a product $\prod_{i=1}^s P(\beta_i)$. Hence under notations of
Formula~\ref{Rem2} $\Discr(f(u,y)-v) =
c\prod_{j=1}^{n-1}[f(u,\gamma_j(u))-v]$ where $c$ is a nonzero
constant and Formula~\ref{Rem1} follows.
\end{proof}

\begin{Formula}\label{Rem3}
Let $f(x,y)=y^N+a_1(x)y^{N-1}+\dots+a_N(x)\in \mathbf C\{x\}[y]$.
Assume that all nonzero roots of the polynomial $f(0,y)$ are simple.
Then the discriminant of the mapping $(x,f):(\bC^2,0)\to(\bC^2,0)$
is given by formula~(\ref{Eq:Rem1}).
\end{Formula}
\begin{proof}
Let $\frac{\partial f}{\partial y}(x,y)=N
\prod_{i=1}^{N-1}(y-\gamma_i(x))$ be the Puiseux factorization of
$y$-partial derivative. Take $\gamma_k(x)$ such that
$\gamma_k(0)\neq 0$. Since $\gamma_k(0)$ is a root of
$\frac{\partial f}{\partial y}(0,y)$ and all nonzero roots of
$f(0,y)$ are simple there is $f(0,\gamma_k(0))\neq 0$.

\medskip

\noindent We get
\begin{eqnarray*}
\Discr(f(u,y)-v) &=&
\mathit{const}\prod_{i=1}^{N-1}[v-f(u,\gamma_i(u))] \\
&=& \mathit{const}\prod_{\gamma_i(0)\neq0}[v-f(u,\gamma_i(u))]\,
\prod_{\gamma_i(0)=0} [v-f(u,\gamma_i(u))] \\
&=& \mathit{unit}\prod_{\gamma_i(0)=0} [v-f(u,\gamma_i(u))]
\end{eqnarray*}
which is, up to a unit, (\ref{Rem2}).
\end{proof}

\section{Jacobian Newton diagrams of irreducible series}
\label{jacobian branch}

\noindent In this section we consider
mappings~$(l,f):(\bC^2,0)\to(\bC^2,0)$ under additional assumption
that $f$~is an irreducible singular power series. Then the curve
$f=0$ is often called a \textit{plane singular branch}.

\noindent Consider
$$
  S(f)=\{(f,g)_0:\;g\in\bC\{x,y\}\mbox{ and $f$ does not divide $g$}\}\;.
$$

\noindent Clearly $0\in S(f)$ (take $g=1$) and if $a,b\in S(f)$ then
$a+b\in S(f)$ since the intersection number is additive, so $S(f)$
is a semigroup, called the {\it semigroup of the branch\/} $f=0$.

\medskip

\noindent For any regular curve $l=0$ the semigroup $S(f)$ has the
$(f,l)_0$-{\em minimal system of generators}
$\bb_0,\bb_1,\dots,\bb_h$ defined by conditions
\begin{itemize}
\item[\rm(i)]   $\bb_0=(f,l)_0$,
\item[\rm(ii)]  $\bb_k=\min(S(f)\setminus (\bN\,\bb_0+\dots+\bN\,\bb_{k-1}))$,
\item[\rm(iii)] $S(f)=\bN\,\bb_0+\dots+\bN\,\bb_h$.
\end{itemize}

\noindent The sequence of generators can be characterized in purely
arithmetical terms. Let us recall (see \cite{Bresinsky},
\cite{Zariski} for the generic case ($\bb_0=\ord f$) and~\cite{GP2}
for the case when the curves $f=0$, $l=0$ are tangent)
\begin{Theorem}\label{Th:Z}
Let $\bb_0,\bb_1,\dots,\bb_h$ be a sequence of positive integers.
Set $n_k=\gcd(\bb_0,\dots,\bb_{k-1})/\gcd(\bb_0,\dots,\bb_k)$ for
$k\in \{1,\dots,h\}$. Then the following two conditions are
equivalent.
\begin{itemize}
\item[\rm(i)] There is a singular branch $f=0$ and a regular curve
$l=0$ such that $\bb_0,\bb_1,\dots,\bb_h$ is the $(f,l)_0$-minimal
system of genera\-tors of the semigroup $S(f)$,
\item[\rm(ii)] the sequence $\bb_0,\bb_1,\dots,\bb_h$ satisfies
the conditions:
\begin{itemize}
\item[\rm(Z$_1$)] $n_k>1$ for $k\in \{1,\dots,h\}$ and $n_1\cdots n_h=\bb_0$,
\item[\rm(Z$_2$)] $n_k\bb_k<\bb_{k+1}$ for $k\in \{1,\dots,h-1\}$.
\end{itemize}
\end{itemize}
\end{Theorem}
\medskip

\noindent Now we can state the result due to~\cite{Smith},
\cite{Merle} and~\cite{Ephraim}.
\begin{Theorem}[Smith--Merle--Ephraim]\label{Th:M}
Suppose that $f=0$ is a singular branch and $l=0$ is a regular curve.
Let $\bb_0,\dots,\bb_h$ be the $(f,l)_0$-minimal system of
genera\-tors of the semigroup $S(f)$. Then with the notation
introduced above
\begin{equation}\label{Eq:M}
  \Nj(l,f)=\sum_{k=1}^h\Teisssr{(n_k-1)\bb_k}
                               {(n_k-1)n_1\cdots n_{k-1}}{18}{9}\;.
\end{equation}
\end{Theorem}

\medskip

\noindent If $\bb_0,\dots,\bb_h$ is the sequence satisfying the
conditions~ (Z$_1$) and~(Z$_2$) of Theorem~\ref{Th:Z} then we will
write ${\cal M}(\bb_0,\dots,\bb_h)$ for the Newton
diagram~(\ref{Eq:M}) and we call it the {\em Merle type diagram}.
Let us note that the Newton diagram in formula~(\ref{Eq:M}) is
written in the canonical form because quotients of the inclinations
of successive elementary Newton diagrams, which are
$\bb_{k+1}/(n_k\bb_k)$, are greater than~$1$ by Theorem~\ref{Th:Z}.

\medskip

\noindent Let us look at Example~\ref{Ex:1} in the light of
Theorem~\ref{Th:M}. The curve $f=0$ has the semigroup $S(f)=\bN
2+\bN 5$. There is $(f,l_1)_0=2$, $(f,l_2)_0=5$, $(f,l_3)_0=4$ and
is easy to verify that $\Nj(l_1,f)={\cal M}(2,5)$, $\Nj(l_2,f)={\cal
M}(5,2)$ and $\Nj(l_3,f)={\cal M}(4,2,5)$.

\begin{Theorem}\label{Th1}
Let $\Delta=\sum_{i=1}^h \Teis{L_i}{M_i}$ be a convenient Newton
diagram written in its canonical form. Put $H_0=1$,
$H_i=1+M_1+\dots+M_i$ for  $i\in \{1,\dots,h\}$ and $C_0=H_h$,
$C_i=H_{i-1}L_i/M_i$ for $i\in \{1,\dots,h\}$.
Then $\Delta$ is a Merle type diagram  if and only if the arithmetic
conditions (i)--(iii) are satisfied
\begin{itemize}
\item[(i)] the quotients  $H_i/H_{i-1}$ are integers for $i\in\{2,\dots,h\}$,
\item[(ii)] the quotients $C_i$ are integers for $i\in \{1,\dots,h\}$,
\item[(iii)] $\gcd(C_0,\dots,C_i)=C_0/H_i$ for $i\in \{1,\dots,h\}$.
\end{itemize}
\noindent Moreover in such a case $\Delta={\cal M}(C_0,\dots,C_h)$.
\end{Theorem}

\noindent \begin{proof} Assume that $\Delta$ is a Merle type diagram
${\cal M}(\bb_0,\dots,\bb_h)$. Then $L_i=(n_i-1)\bb_i$ and
$M_i=(n_i-1)n_1\cdots n_{i-1}$ for $i\in \{1,\dots,h\}$. We have the
equality $H_i=n_1\cdots n_i$. Indeed $H_1=1+M_1=1+(n_1-1)=n_1$ and
$H_{i+1}=H_i+M_{i+1}=n_1\cdots n_i+(n_{i+1}-1)n_1\cdots
n_i=n_1\cdots n_{i+1}$ by the inductive hypothesis. It follows that
$H_i/H_{i-1}=n_i$ hence condition~(i) is satisfied.

\medskip

\noindent It also follows that
$C_i=\frac{H_{i-1}L_i}{M_i}=\frac{n_1\cdots
n_{i-1}(n_i-1)\bb_i}{(n_i-1)n_1\cdots n_{i-1}}=\bb_i$. Hence
condition (ii) is also satisfied.

\medskip

\noindent It follows directly from the definition of the sequence
$n_i$ that $\gcd(\bb_0,\dots,\bb_i)=\bb_0/(n_1\cdots n_i)$ for
$1\leq i\leq h$. Moreover by condition~ (Z$_1$ of Theorem~\ref{Th:Z}
there is $n_1\cdots n_h=\bb_0$ which gives $C_0=H_h=\bb_0$. Thus  by
condition~ (Z$_1$) of Theorem~\ref{Th:Z}
$\gcd(C_0,\dots,C_i)=C_0/H_i$ for $i\in\{1,\dots,h\}$.

\medskip

\noindent
Now assume that conditions (i)--(iii) hold true for the Newton diagram $\Delta$.
We will show that the sequence $C_0,\dots,C_h$ satisfies arithmetical conditions of
Theorem~\ref{Th:Z}. It follows from~(iii) that $n_i:=
\frac{\gcd(C_0,\ldots,C_{i-1})}{\gcd(C_0,\ldots,C_{i})}=H_i/H_{i-1}$
for $i\in \{1,\dots,h\}$. Thus $n_i>1$ for $i\in \{1,\dots,h\}$ and
$n_1\cdots n_h=C_0$.

\medskip

\noindent Because $\Delta$ is written in canonical form there is
$L_i/M_i<L_{i+1}/M_{i+1}$ for $i\in \{1,\dots,h-1\}$. Multiplying
these inequalities by $n_i H_{i-1}=H_{i}$ we get $n_iH_{i-1}
L_i/M_i<H_{i}L_{i+1}/M_{i+1}$ which is equivalent with
$n_iC_i<C_{i+1}$ for $i\in \{1,\dots,h-1\}$.

\medskip

\noindent We checked that the sequence $C_0,\dots,C_h$ satisfies
conditions~(Z$_1$) and~(Z$_2$) of Theorem~\ref{Th:Z}. Moreover
looking at the first part of the proof it is easy to see that
$\Delta={\cal M}(C_0,\dots,C_h)$.
\end{proof}
\medskip


\section{Discriminant criterion of irreducibility}
\begin{Theorem}\label{Th2}
Let $f=0$ be a plane singular curve and let $l=0$ be a regular
curve.
Then $f$ is irreducible if and only if $\Nj(l,f)$ is a Merle type
diagram. Moreover if $\Nj(l,f)={\cal M}(\bb_0,\dots,\bb_h)$ then
$f=0$ has the semigroup $S(f)=\bN\bb_0+\cdots +\bN\bb_h$.
\end{Theorem}

\begin{Example}\label{Ex1} Let $l=x$ and $f=y^n-x^m$. Then by~(\ref{Eq:Rem2})
$D(u,v)=(v+u^m)^{n-1}$, hence
$\Nj(x,f)=\Teisssr{(n-1)m}{n-1}{10}{5}$. Under notation of
Theorem~\ref{Th1} there is $C_0=H_1=n$, $C_1=m$ and conditions (i)
and (ii) of Theorem \ref{Th1} are clearly satisfied. Condition~(iii)
reduces to $\gcd(m,n)=1$ and it is well-known that the curve
$y^n-x^m=0$ is irreducible if and only if $m$ and $n$ are co-prime.
\end{Example}

\noindent The following two examples are taken from \cite{Kuo} (see
also \cite{Abhyankar}).

\begin{Example}\label{Ex2}
Let $f=(y^2-x^3)^2-x^7$. Then
$\jac(x,f)=4y(y^2-x^3)=4y(y-x^{3/2})(y+x^{3/2})$. By Formula
\ref{Rem2} we get $D(u,v)=(v-u^6+u^7)(v+u^7)^2$. Hence
$\Nj(x,f)=\Teis{6}{1}+\Teis{14}{2}$. Under notation of
Theorem~\ref{Th1} there is $H_1=1+1=2$, $C_0=H_2=1+1+2=4$,
$C_1=6/1=6$, $C_2=H_1\cdot 14/2=14$ and because
$\gcd(C_0,C_1,C_2)=2\neq1$, it follows that $\Nj(x,f)$ is not a
Merle type diagram. Therefore $f$ is not irreducible.
\end{Example}

\begin{Example}\label{Ex3}
Let $f(x,y)=(y^2-x^3)^2-x^5y$. By Formula \ref{Rem1}
$D(u,v)=-256v^3+256u^6v^2+288u^{13}v-256u^{19}-27u^{20}$ (we
computed the discriminant using {Sage}) and the Newton diagram of
the discriminant is $\Nj(x,f)=\Teis{6}{1}+\Teis{13}{2}$. It is easy
to check that $\Nj(x,f)$ is a Merle type diagram ${\cal M}(4,6,13)$.
Therefore $f$ is irreducible with semigroup $S(f)=\bN4+\bN6+\bN13$.
\end{Example}

\begin{Example}\label{Ex4}
Let $f(x,y)=x^8+(x^2+y^3)^3$. The jacobian Newton diagram  of $(x,f)$
is $\Nj(x,f)=\Teis{12}{2}+ \Teis{48}{6}$ which is not
a Merle type diagram. Note that in this example $x=0$ is not transverse to
$f(x,y)=0$.
\end{Example}

\begin{Corollary}
\label{coro} Let $f(x,y)=y^N+a_1(x)y^{N-1}+\dots+a_N(x)\in \mathbf
C[x,y]$. Assume that the curve $f(x,y)=0$ intersects $x=0$ only at
the point $(0,y_0)$. Then the curve $f(x,y)=0$ is analytically
irreducible at $(0,y_0)$ if and only if the Newton diagram of
$\Discr(f(u,y)-v)$ is a Merle type diagram.
\end{Corollary}

\noindent \begin{proof} Put $\tilde{f}(x,y)=f(x,y+y_0)$. Then
$f(x,y)=0$ is analytically irreducible at $(0,y_0)$ if and only if
$\tilde{f}(x,y)=0$ is analytically irreducible at $(0,0)$. Since
$\Discr(f(u,y)-v)=\Discr(\tilde{f}(u,y)-v)$ the result follows from
Formula~\ref{Rem3}.
\end{proof}

\section{Proof of Theorem \ref{Th2}}
The proof is based on Theorem 1 of \cite{GG}:
\begin{Theorem}\label{Th:GB}
Let $f,g\in\bC\{x,y\}$ be such that $\Nj(x,f)=\Nj(x,g)$. Assume that
$x=0$ is transverse to the curves $f=0$ and $g=0$. If $f$ is
irreducible then $g$ is also irreducible.
\end{Theorem}
and on the following lemma
\begin{Lemma}\label{L1}
Let $f$ where $f(0,y)=y^n+\mbox{higher order terms }$ be a
convergent power series and let $N$ be a positive integer. Put
$\tilde f(x,y)=f(x^N,y)$. Then
\begin{itemize}
\item[(i)]  if $N$ and $n$ are coprime integers
            then $f$ is irreducible if and only if $\tilde f$ is irreducible,
\item[(ii)] if $N>n$ then $\tilde f=0$ is transverse to $x=0$,
\item[(iii)] $\Nj(x,\tilde f)=L(\Nj(x,f))$ where $L:\bR^2\to\bR^2$
            is a linear automorphism given by
            $L(i,j)=(Ni,j)$.
\end{itemize}
\end{Lemma}

\noindent \begin{proof}

\noindent \textit{Proof of (i).} Assume that $f=f_1f_2$. Then
$\tilde f(x,y)=f_1(x^N,y)f_2(x^N,y)$. It follows that if $\tilde
f$ is irreducible then $f$ is irreducible.

\medskip

\noindent In order to show the implication in opposite direction
assume that $f$ is irreducible. Recall (see Theorem 2.1, Chapter IV,
\cite{Walker}) that the curve $f=0$ where $\ord f(0,y)=n$ is a
branch if and only if there exists a convergent power series
$\phi(t)$ such that $f(t^n,\phi(t))=0$ and the greatest common
divisor of the set $\{n\}\cup\supp\phi$ equals 1.

\medskip

\noindent Let $\phi(t)$ be such a series and let
$\tilde\phi(t)=\phi(t^N)$. Then $\tilde
f(t^n,\tilde\phi(t))=f(t^{nN},\phi(t^N))=0$ and since $n$ and $N$
are co-prime the greatest common divisor of the set
$\{n\}\cup\supp\tilde\phi=\{n\}\cup N\cdot\supp\phi$ equals 1.
Consequently $\tilde f=0$ is a branch.

\medskip \noindent \textit{Proof of (ii).} By the assumption $N>n$
the homogeneous initial part of the series $f(x^N,y)$ is $y^n$. This
gives (ii).

\medskip \noindent \textit{Proof of (iii).} Let $\frac{\partial
f}{\partial y}(x,y) = \mathit{unit}\prod_{j=1}^{n-1}[y-\gamma_j(x)]$
be the Newton-Puiseux factorization of $\frac{\partial f}{\partial
y}$. By Formula \ref{Rem2} the discriminant  of the mapping $(x,f)$
has an equation $D(u,v)=\prod_{j=1}^{n-1}[v-f(u,\gamma_j(u))]$.
Because $\frac{\partial \tilde f}{\partial y}(x,y)=\frac{\partial
f}{\partial y}(x^N,y)$ there is $\frac{\partial \tilde f}{\partial
y}(x,y) = \mathit{unit}\prod_{j=1}^{n-1}[y-\gamma_j(x^N)]$ and
consequently the discriminant of the mapping $(x,\tilde f)$ has an
equation $\tilde D(u,v)=\prod_{j=1}^{n-1}[v-\tilde
f(u,\gamma_j(u^N))] =
\prod_{j=1}^{n-1}[v-f(u^N,\gamma_j(u^N))]=D(u^N,v)$. Comparing
$\Delta_D$ with $\Delta_{\tilde D}$ we get (iii).
\end{proof}

\bigskip

\noindent Now let us prove Theorem~\ref{Th2}. Suppose that
$\Nj(l,f)=\Nj(l,g)$ where $g$ is an irreducible power series.
Applying an analytic change of coordinates we may assume that $l=x$.
Take an integer $N>0$ such that conclusions of~(i) and~(ii) of
Lemma~\ref{L1} are satisfied for $\tilde f(x,y)=f(x^N,y)$ and
$\tilde g(x,y)=g(x^N,y)$. It follows from~(iii) of Lemma~\ref{L1} that
$\Nj(x,\tilde f)=\Nj(x,\tilde g)$.
Since $\tilde f$ and $\tilde g$ satisfy assumptions of
Theorem~\ref{Th:GB}, $\tilde f$ is an irreducible power series. Hence by
(i) of Lemma~\ref{L1} $f$ is also irreducible.

\section{Discriminant criterion of irreducibility at infinity}

\noindent Let $p(x,y)$ be a complex polynomial of degree $n>0$. Let
$C\subset \mathbf P^2(\mathbf C)$ be the projective closure of the
curve $p(x,y)=0$.
 Assume that $C$ intersects the line at infinity
at only one point $Q$. The purpose of this section is to give a
criterion for local analytical irreducibility of the curve $C$ at
$Q$ without passing to local coordinates centered at $Q$. For this
we need some notions.

\medskip

\noindent Let $g(x,y)$ be a polynomial of positive degree such that
$g(x,0) \neq 0$ and $g(0,y)\neq 0$ (in other words its Newton
diagram is convenient). Let ${\cal P}_0(g)$ be the boundary in
$\mathbf R^2_{+}$ of $\Delta_g$ and ${\cal P}_{\infty}(g)$ be the
boundary in $\mathbf R^2_+$ of $\Delta_{\infty}(g) = \mbox{Convex
Hull}\;(\supp g \, \cup \, \{(0,0)\})$. We call these sets the {\em
Newton polygon of $g$ at zero} and the {\em Newton polygon of $g$ at
infinity} respectively.

\medskip
\noindent

\begin{Theorem}
\label{infty} Let $p(x,y)$ be a complex polynomial of degree $n>0$
without multiple factors and let $C$ be the projective closure of
$p(x,y)=0$. Assume that $C$ intersects the line at infinity at only
one point $Q\neq (0:1:0)$. Put $D_{\infty}(x,t):=\Discr(p(x,y)-t)$
and let $L:\mathbf Z^2 \longrightarrow \mathbf Z^2$ be the affine
transformation defined by $L(i,k)=(n(n-1)-i-nk,k)$. Then the curve
$C$ is analytically irreducible at $Q$ if and only if $L({\cal
P}_{\infty}(D_{\infty}))$ is the Newton polygon at zero of a Merle
type diagram.
\end{Theorem}

\noindent \begin{proof}
Let $P(x,y,z)=z^np\left(\frac{x}{z},\frac{y}{z}\right)$ be a
homogeneous equation of the curve $C$. Assume that $C$
intersects the line at infinity only at $Q=(1:y_0:0)$.
Then $p(x,y)=P(x,y,1)$ and $f(y,z):=P(1,y,z)=0$ is the affine equation of $C$
in coordinates $y,z$. In these coordinates the point $Q$ becomes $(y_0,0)$.
Since the curve $C$ intersects the line $z=0$ only at $Q$ the polynomial
$f(y,z)$ satisfies the assumptions of Corollary \ref{coro}.

\medskip

\noindent Put $D(x,z,t)=\mathrm{Disc}_y(P(x,y,z)-t)$. We have that
$D_{\infty}=D(x,1,t)$ and
$D_{0}:=\mathrm{Disc}_y(f(z,y)-t)=D(1,z,t)$. Since $P(x,y,z)$ is a
homogeneous polynomial of degree $n$, giving to the variable $t$ the
weight $n$, and the other variables the weight $1$, the polynomial
$D(x,z,t)$ is quasi-homogeneous of degree $n(n-1)$ (see Theorem
10.9, Chapter I, \cite{Walker}). In particular any term
$c_{ijk}x^iz^jt^k$ of $D(x,z,t)$ corresponds with the point
$(i,j,k)$ of the hyperplane $\Pi\equiv i+j+nk=n(n-1)$. Moreover such
point determines the term $c_{ijk}x^it^k$ of $D_{\infty}$ and the
term $c_{ijk}z^jt^k$ of $D_{0}$. Put $L:\mathbf Z^2\longrightarrow
\mathbf Z^2$ defined by $L(i,k)=(n(n-1)-i-nk,k)$.
 Then
$\supp D_{0} = L(\supp D_{\infty})$. The last relation gives ${\cal
P}_{0}(D_0)=L({\cal P}_{\infty} (D_{\infty}))$. By Corollary
\ref{coro} the curve $C$ is analytically irreducible at $Q$ if and
only if the Newton diagram $\Delta_{D_0}$ is a Merle type
diagram.\end{proof}

\medskip

\begin{Remark}
Let $p(x,y)$ be a polynomial of degree $n$ which has one point at
infinity different from $(0:1:0)$. Let us denote by $q$ the maximal
inclination of ${\cal P}_{0}(D_{0})$. After \cite{Ploski}, the
Abhyankar-Moh inequality (see \cite{Abhyankar-Moh}) is equivalent to
$q<n$. Note also that the Abhyankar-Moh inequality is equivalent to
equisingularity at infinity of the family $p(x,y)-t=0$. After
\cite{Krasinski} this is also equivalent to the statement that all
segments of ${\cal P}_{\infty}(D_{\infty})$ have positive slopes.
\end{Remark}

\medskip

\begin{Example}
Let $p(x,y)=x+(x+y^3)^3$ be a polynomial in $\mathbf C[x,y]$ which
corresponds to the  projective curve $C$ defined by
$P(x,y,z)=xz^8+(xz^2+y^3)^3 =0$. The only point at infinity of $C$
is $Q=(1:0:0)$. Moreover $D_{\infty}=(x+x^3-t)^2(x-t)^6$ and ${\cal
P}_{\infty}(D_{\infty})$ has only two segments joining the point
$(0,8)$ to $(6,6)$ and this one to $(12,0)$. The transformation of
${\cal P}_{\infty}(D_{\infty})$ by $L(i,k)=(72-i-9k,k)$ is a polygon
of two segments joining the point $(0,8)$ to $(12,6)$ and this one
to $(60,0)$. This polygon is the Newton polygon of
$\Delta=\Teis{12}{2}+ \Teis{48}{6}$. Since $\Delta$ is not a Merle
type diagram, by Theorem \ref{infty} the curve  $C$ is not
analytically irreducible at $Q$. Observe that the local equation
$P(1,y,z)=0$ of the curve $C$ was studied in Example \ref{Ex4}.
\end{Example}

\section*{Appendix} The purpose of this section is the proof of Formula \ref{Rem2}.
We establish a more general statement. Using basic properties of
direct image we show that if $F=(x,f)$ and
$h(x,y)=\prod_{j=1}^{m}[y-\gamma_j(x)]$ is the Newton-Puiseux
factorization of the power series $h$ then the direct image
$F_{*}(h=0)$ has an equation
\begin{equation}\label{Eq:Ap}
\prod_{j=1}^{m}[v-f(u,\gamma_j(u))]=0.
\end{equation}

\begin{Lemma}\label{L1:Ap}
Let $\phi(t)\in\bC\{t\}$, $\phi(0)=0$ be a convergent power series
and let $m$ be a positive integer. Let $V\subset (\bC^2,0)$ be the
image of the mapping $(\bC,0)\ni t\to(t^m,\phi(t))\in(\bC^2,0)$.
Then $V$ has an analytic (not necessarily reduced) equation $g=0$
where
$$ g(x,y)=\prod_{\epsilon^m=1}[y-\phi(\epsilon x^{1/m})].
$$
\end{Lemma}

\noindent Lemma~\ref{L1:Ap} is an easy corollary of:

\begin{Theorem}[Puiseux Theorem]
Let $h(x,y)$ be an irreducible power series and let $n=\ord
f(0,y)<\infty$. Then there exists $\phi(t)\in\bC\{t\}$, $\phi(0)=0$
such that $h(t^n,\phi(t))=0$. Moreover $h$, up to an invertible
factor, is equal to
\begin{equation}\label{Puiseux}
H(x,y):=\prod_{\epsilon^n=1}[y-\phi(\epsilon x^{1/n})] .
\end{equation}
Conversely, if $(\bC,0)\ni t\to(t^n,\phi(t))\in(\bC^2,0)$
is an analytic parametrization of a plane branch then $H(x,y)=0$ is
its reduced equation.
\end{Theorem}

\noindent \begin{proof}(of Lemma \ref{L1:Ap}) Let $d$ be the
greatest common divisor of $\{m\}\cup\,\supp\phi$. Then there exists
a power series $\phi_0$ such that $\phi(t)=\phi_0(t^d)$.

\noindent Put $n=m/d$. Since the greatest common divisor of
$\{n\}\cup\,\supp\phi_0$ equals 1 the mapping $t\to (t^n,\phi_0(t))$
is an analytic parametrization of $V$. Let $\epsilon_m$ be the
$m$-th primitive root of unity and let $\epsilon_n=\epsilon_m^d$. By
Puiseux' Theorem $V$~has an equation
$h(x,y)=\prod_{i=1}^n[y-\phi_0(\epsilon_n^i x^{1/n})]=0$.

\noindent Put $g=h^d$. Then
$g(x,y)= \prod_{i=1}^n[y-\phi_0(\epsilon_n^i x^{1/n})]^d =
         \prod_{i=1}^{nd}[y-\phi_0(\epsilon_n^i x^{1/n})] =
         \prod_{i=1}^{m}[y-\phi_0((\epsilon_m^i x^{1/m})^d)] =
         \prod_{i=1}^{m}[y-\phi(\epsilon_m^i x^{1/m})]$ and we get
Lemma \ref{L1:Ap}. Observe that $g$ is not reduced in the
case~$d>1$.
\end{proof}

\medskip

\noindent Now we prove~(\ref{Eq:Ap}) under assumption that $h=0$ is a branch.

\begin{Lemma}\label{Ap1}
Let $h(x,y)=\prod_{i=1}^{n}[y-\gamma_i(x)]$ be an irreducible power
series. Then the direct image of the branch $h=0$ by the mapping
$(x,f):(\bC^2,0)\to(\bC^2,0)$ has an equation
$\prod_{i=1}^{n}[v-f(u,\gamma_i(u))]=0$.
\end{Lemma}

\noindent \begin{proof}
By Puiseux' Theorem we may assume that
$\gamma_i(x)=\phi(\epsilon_n^i x^{1/n})$ for $i=1,\dots,n$
where $t\to (t^n,\phi(t))$ is an analytic
parametrization of a branch $h=0$ and $\epsilon_n$
is the $n$-th primitive root of unity.

\noindent Thus the image of the curve $h=0$ by the mapping $(x,f)$
is given by $(u,v)=(t^n,f(t^n,\phi(t)))$. By Lemma~\ref{L1:Ap} it
has an equation $g(u,v)=\prod_{i=1}^n[v-f(u,\phi(\epsilon_n^i
u^{1/n}))]=\prod_{i=1}^{n}[v-f(u,\gamma_i(u))]=0$. Moreover
$(g,u)_0=(h,x)_0=n$ which shows, after Projection Formula (see page
64 of \cite{Casas}), that $g=0$ is the direct image of the curve
$h=0$.  \end{proof}

\medskip

\noindent
Because the equation of the direct image of a curve
$h_1h_2=0$ is the product of equations of direct images of $h_i=0$
($i=1,2$) formula~(\ref{Eq:Ap}) holds also in the case when $h$ is a
reducible power series.

\medskip
\noindent
{\small Evelia Rosa Garc\'{\i}a Barroso\\
Departamento de Matem\'atica Fundamental\\
Facultad de Matem\'aticas, Universidad de La Laguna\\
38271 La Laguna, Tenerife, Espa\~na\\
e-mail: ergarcia@ull.es}

\medskip

\noindent {\small Janusz Gwo\'zdziewicz\\
Department of Mathematics\\
Technical University \\
Al. 1000 L PP7\\
25-314 Kielce, Poland\\
e-mail: matjg@tu.kielce.pl}

\end{document}